# On Bergman type spaces of holomorphic functions and the density, in these spaces, of certain classes of singular functions

by

**T. Hatziafratis, K. Kioulafa, and V. Nestoridis**

### Abstract

We consider Bergman spaces and variations of them on domains $\Omega$ in one or several complex variables. For certain domains $\Omega$ we show that the generic function in these spaces is totaly unbounded in $\Omega$ and hence non-extendable. We also show that generically these functions do not belong – not even locally – in Bergman spaces of higher order. Finally, in certain domains $\Omega$, we give examples of bounded non-extendable holomorphic functions – a generic result in the spaces $A^s(\Omega)$ of holomorphic functions in $\Omega$ whose derivatives of order $\leq s$ extend continuously to $\overline{\Omega}$ ($0 \leq s \leq \infty$).

**1. Introduction.** An important problem in complex analysis is whether there exists a holomorphic function $f$, in a given open set $\Omega$ in $\mathbb{C}^n$, which is singular at every boundary point of $\Omega$ in the sense that whenever $U$ and $V$ are open subsets of $\mathbb{C}^n$, with $U$ being connected and $\varnothing \neq V \subseteq U \cap \Omega \neq U$, then there is no holomorphic function $F$ in $U$ which extends $f|_V$, i.e., $F(z) = f(z)$ for $z \in V$. See for example [1], [3], [5], [6], [8], [9] and [11]. Also the problem of constructing singular functions with specific properties – for example satisfying certain growth conditions near the boundary or having certain smoothness upto the boundary – has been studied in various directions. See for example [7], [8], [10], and [11].

In this paper we will show – under certain restrictions on the open set $\Omega$ – that the set of the $\mathcal{O}L^p$ functions (holomorphic and $L^p$ with respect to Lebesgue measure) in $\Omega$, which are totally unbounded, is dense and $\mathcal{G}_\delta$ in the space $\mathcal{O}L^p(\Omega)$. In fact we work mostly with the space $\bigcap\limits_{p<q} \mathcal{O}L^p(\Omega)$ endowed with its natural topology. For example we show that for some open sets $\Omega$, the set of the functions in $\bigcap\limits_{p<q} \mathcal{O}L^p(\Omega)$, which are not in $\mathcal{O}L^q(B(\zeta,\varepsilon) \cap \Omega)$ for any $\zeta \in \partial\Omega$ and any $\varepsilon > 0$, is dense and $\mathcal{G}_\delta$ in the space $\bigcap\limits_{p<q} \mathcal{O}L^p(\Omega)$. We mention that the papers [2] and [10] contain some results which are related to questions studied here. In fact it was the paper [2] which gave us the idea to pose these questions.

On the other hand, as it is well-known, a singular function in $\Omega$ can be $C^\infty$ upto the boundary of $\Omega$. In this direction, we will use a theorem from [12], to show that for some pseudoconvex open sets $\Omega$, the set of the functions in $A^s(\Omega)$ (holomorphic in $\Omega$, whose derivatives of order $\leq s$ extend continuously to $\overline{\Omega}$, for a fixed $s \in \{0,1,2,...\} \cup \{\infty\}$), which are singular at every boundary point of $\Omega$, is dense and $\mathcal{G}_\delta$ in the space $\mathrm{A}^s(\Omega)$ (in the natural topology of this space).

**2. Preliminaries.** The main results of this paper are Theorems 4 and 6 below, and some applications which are Theorems 10, 12, 14 and 16. To prove Theorem 4, we will use the following result from [15].

**Theorem (i).** *Let $\mathcal{V}$ be a topological vector space over $\mathbb{C}$, $X$ a non-empty set, and let $\mathbb{C}^X$ denote the vector space of all complex-valued functions on $X$. Suppose $T : \mathcal{V} \to \mathbb{C}^X$ is a linear operator with the property that, for every $x \in X$, the functional $T_x : \mathcal{V} \to \mathbb{C}$, defined by $T_x(f) := T(f)(x)$, for $f \in \mathcal{V}$, is continuous. Let $S = \{f \in \mathcal{V} : T(f) \text{ is unbounded on } X\}$. Then either $S = \varnothing$ or $S$ is dense and $\mathcal{G}_\delta$ set in the space $\mathcal{V}$.*

Let us also observe that the above theorem holds under the weaker assumption that the operator $T$ is only sublinear (not necessarily linear) in the following sense: If $f, g \in \mathcal{V}$ and $\lambda \in \mathbb{C}$ then

$$|T(f+g)| \le |T(f)| + |T(g)| \text{ and } |T(\lambda f)| = |\lambda| |T(f)|.$$

Indeed, if $S \ne \varnothing$, we will show that $S$ is dense in $\mathcal{V}$. Let $\varphi \in S$. For every $f \in \mathcal{V} - S$, the sequence $f + \frac{1}{k}\varphi$, $k = 1,2,3,...$, converges to $f$, in $\mathcal{V}$. Notice that the function $T(f)$ is bounded (since $f \in \mathcal{V} - S$) and $T(\varphi)$ is unbounded. We claim that (for each fixed $k$) the function $T\left(f + \frac{1}{k}\varphi\right)$ is unbounded. For if it were bounded, the sublinearity of $T$,

$$\left|\frac{1}{k}T(\varphi)\right| = \left|T\left(\frac{1}{k}\varphi\right)\right| \le \left|T\left(f + \frac{1}{k}\varphi\right)\right| + |T(-f)| = \left|T\left(f + \frac{1}{k}\varphi\right)\right| + |T(f)|,$$

would imply that $T(\varphi)$ is bounded. Thus $T\left(f + \frac{1}{k}\varphi\right)$ is unbounded, i.e., $f + \frac{1}{k}\varphi \in S$, and the density of $S$ in $\mathcal{V}$ follows. The assertion that $S$ is $\mathcal{G}_\delta$ in $\mathcal{V}$ can be proved as in the case $T$ is linear, using the continuity of the functionals $T_x$, $x \in X$.

We point out that in the previous Theorem (i), Baire's theorem is not used and the space $\mathcal{V}$ is not assumed to be complete.

We will need the above observation in the proof of Theorem 6. We will also use the following theorem which is proved in [12]. See also [8] for a related result.

**Theorem (ii).** *Let $\Omega \subset \mathbb{C}^n$ be an open set and let $\mathcal{V}$ be a vector subspace of $\mathcal{O}(\Omega)$. Suppose that in $\mathcal{V}$ there is defined a complete metric whose topology makes $\mathcal{V}$ a topological vector space and such that convergence in $\mathcal{V}$ implies pointwise convergence in $\mathcal{O}(\Omega)$. If for every pair of balls $(B,b)$ with $b \subset\subset B \cap \Omega \ne B$, there exists $f_{(B,b)} \in \mathcal{V}$ such that the restriction $f_{(B,b)}\big|_b$ (of the function $f_{(B,b)}$ to $b$) does not have any bounded holomorphic extention to $B$, then the set of the functions $g \in \mathcal{V}$ which are singular at every boundary point of $\Omega$ is dense and $\mathcal{G}_\delta$ in $\mathcal{V}$.*

**The spaces $\mathcal{O}L^p(\Omega)$.** Let $\Omega \subset \mathbb{C}^n$ be a bounded open set. We recall that for $p \ge 1$, the Bergman space $\mathcal{O}L^p(\Omega)$ is defined to be the set of holomorphic functions $f : \Omega \to \mathbb{C}$ such that

$$\|f\|_p := \left(\int_\Omega |f(z)|^p \, dv(z)\right)^{1/p} < +\infty,$$

where $dv$ is the Lebesque measure in $\mathbb{C}^n$. Then the quantity $\|\cdot\|_p$ is a norm, and with this norm, $\mathcal{O}L^p(\Omega)$ is a Banach space. Indeed, recalling the inequality

$$\sup_{K}\left|\frac{\partial^{\alpha} f}{\partial z^{\alpha}}\right| \leq c(\alpha, K)\|f\|_{1}, \text{ for } f \in \mathcal{O}L^{1}(\Omega) \text{ (see [6])},$$

where $K$ is a compact subset of $\Omega$, and $c(\alpha, K)$ is a constant depending on $K$ and the multi-index $\alpha$, it follows that if a sequence $f_{k} \in \mathcal{O}L^{1}(\Omega)$ converges to $f$, in the $L^{1}(\Omega)$−norm, then the convergence is uniform on compact subsets of $\Omega$. In particular, $\mathcal{O}L^{1}(\Omega)$ is closed subspace of $L^{1}(\Omega)$, and, more generally, $\mathcal{O}L^{p}(\Omega)$ is closed subspace of $L^{p}(\Omega)$, for $p \geq 1$. Since we assume $\Omega$ to be bounded, $\mathcal{O}L^{q}(\Omega) \subset \mathcal{O}L^{p}(\Omega)$ when $q > p$. Similarly we define the space $\mathcal{O}L^{\infty}(\Omega)$, of bounded holomorphic functions $f : \Omega \to \mathbb{C}$, which becomes a Banach space with the norm $\|f\|_{\infty} = \sup_{z \in \Omega}|f(z)|$.

For a fixed $q > 1$, we will also consider the spaces

$$\bigcap_{p < q} \mathcal{O}L^{p}(\Omega)$$

endowed with the metric

$$d(f, g) := \sum_{j=1}^{\infty} \frac{1}{2^{j}} \frac{\|f - g\|_{p_{j}}}{1 + \|f - g\|_{p_{j}}}, \ f, g \in \bigcap_{p < q} \mathcal{O}L^{p}(\Omega),$$

where $p_{j}$ is a sequence with $1 < p_{1} < p_{2} < \cdots < p_{j} < \cdots < q$ and $p_{j} \to q$ (as $j \to \infty$). Then $\bigcap_{p < q} \mathcal{O}L^{p}(\Omega)$ becomes a complete metric space, its topology being independant of the choice of the sequence $p_{j}$. In fact, a sequence $f_{k}$ converges to $f$, in the space $\bigcap_{p < q} \mathcal{O}L^{p}(\Omega)$, if and only if $\|f_{k} - f\|_{p} \to 0$ for every $p < q$. Thus Baire's theorem hold in $\bigcap_{p < q} \mathcal{O}L^{p}(\Omega)$: A countable intersection of open and dense subsets of $\bigcap_{p < q} \mathcal{O}L^{p}(\Omega)$ is dense and $\mathcal{G}_{\delta}$ in this space. Moreover we point out that the space $\bigcap_{p < q} \mathcal{O}L^{p}(\Omega)$, with the above topology, is also a topological vector space. In particular, if $f_{k}, f \in \bigcap_{p < q} \mathcal{O}L^{p}(\Omega)$ with $d(f_{k}, f) \to 0$ $(k \to \infty)$, and $\lambda_{k}, \lambda \in \mathbb{C}$ with $\lambda_{k} \to \lambda$, then $d(\lambda_{k} f_{k}, \lambda f) \to 0$. Finally we observe that all the above hold in the case $q = \infty$ too, defining the space $\bigcap_{p < \infty} \mathcal{O}L^{p}(\Omega)$, and that this space contains the space of bounded holomorphic functions in $\Omega$:

$$\bigcap_{p < \infty} \mathcal{O}L^{p}(\Omega) \supset \mathcal{O}L^{\infty}(\Omega).$$

**3. Totally unbounded holomorphic functions.** Let $\Omega \subset \mathbb{C}^{n}$ be an open set. We will say that a holomorphic function $f : \Omega \to \mathbb{C}$ is *totally unbounded* in $\Omega$, if for every $w \in \partial\Omega$, every $\delta > 0$, and every connected component $E$ of the set

$$B(w, \delta) \cap \Omega = \{z \in \Omega : |z - w| < \delta\},$$

the function $f|_{E}$ is unbounded, i.e., $\sup_{z \in E}|f(z)| = \infty$. Notice that such a function is *singular* at every point of $\partial\Omega$. More precisely the following proposition holds.

**Proposition.** *Let $\Omega \subset \mathbb{C}^n$ be a an open set and let $f : \Omega \to \mathbb{C}$ be a totally unbounded holomorphic function. Then for every open sets $U, V \subset \mathbb{C}^n$, with $U$ being connected and $\varnothing \neq V \subseteq U \cap \Omega \neq U$, there does not exist a holomorphic function $F$ on $U$ which extends $f\big|_V$, i.e., $F\big|_V = f\big|_V$.*

**Proof.** Suppose – to reach a contradiction – that for some pair of sets $U$ and $V$, there exists a function $F$, which extends $f$ in the way described above. Let $E_1$ be the connected component of $U \cap \Omega$ which contains $V$. Then $F\big|_{E_1} = f\big|_{E_1}$ and $\overline{E_1} \cap \partial\Omega \cap U \neq \varnothing$, so that we may take a point $w \in \overline{E_1} \cap \partial\Omega \cap U$, and a ball $B(w, \delta)$ with $\overline{B(w, \delta)} \subset U$. Then $B(w, \delta) \cap E_1 \neq \varnothing$, and if $c \in B(w, \delta) \cap E_1$ then for the connected component $E$ of the set $B(w, \delta) \cap \Omega$, which contains the point $c$, we have $\sup_{z \in E} |f(z)| = \infty$ (since $f$ is assumed to be totally unbounded). But this contradicts the equation $F\big|_E = f\big|_E$, which follows from the principle of unique analytic continuation, applied to the connected open set $E$ and the fact that open set $E \cap E_1 \neq \varnothing$. This completes the proof.

**Remark.** In the above proof we used that $\overline{E_1} \cap \partial\Omega \cap U \neq \varnothing$. To justify this elementary topological fact, let us observe that, since $U \cap \Omega \neq \varnothing$, $U \cap (\mathbb{C} - \Omega) \neq \varnothing$ and $U$ is connected, it follows that $U \cap \partial\Omega \neq \varnothing$. Let $a \in V$ and $b \in U \cap \partial\Omega$, and let $\Gamma$ be a curve which lies in $U$ and connects the points $a$ and $b$. If $C$ is the connected component of $U \cap \Omega$ which contains $a$, then $C$ is open, $a \in C \cap \Gamma$ and $b \notin C \cap \Gamma$. Since the set $\Gamma$ is connected, we must have $\Gamma \cap \partial C \neq \varnothing$. Then for a point $\tau \in \Gamma \cap \partial C$, we will have $\tau \in \partial\Omega$ and $\tau \in \overline{C}$, and therefore $\overline{C} \cap \partial\Omega \cap U \neq \varnothing$. Finally, since $E_1 \supseteq C$, we obtain that, indeed, $\overline{E_1} \cap \partial\Omega \cap U \neq \varnothing$.

We will show that under certain assumptions on $\Omega$, the set of the functions in the space $\bigcap_{p<q} \mathcal{O}L^p(\Omega)$, which is totally unbounded in $\Omega$, is dense and $\mathcal{G}_\delta$ (in this space). We will also give examples of specific domains in which this $\mathcal{G}_\delta$ – density conclusion holds.

**4. Theorem.** *Let $\Omega \subset \mathbb{C}^n$ be a bounded open set and $q \in \mathbb{R}$, $q > 1$. Suppose that for every point $\zeta \in \partial\Omega$, there exists a function $f_\zeta$ such that*

$$f_\zeta \in \bigcap_{p<q} \mathcal{O}L^p(\Omega) \quad and \quad \lim_{\substack{z \to \zeta \\ z \in \Omega}} f_\zeta(z) = \infty.$$

*Then the set of the functions $g$ in the space $\bigcap_{p<q} \mathcal{O}L^p(\Omega)$, which are totally unbounded in $\Omega$, is dense and $\mathcal{G}_\delta$ in this space. Also the set of the functions $h$ in the space $\bigcap_{p<q} \mathcal{O}L^p(\Omega)$, which are singular at every boundary point of $\Omega$ is dense and $\mathcal{G}_\delta$ in this space.*

**Proof.** Let us fix a pair $(B, b)$, where $B$ is a «small» open ball whose center lies on $\partial\Omega$ and $b$ is a «smaller» open ball with $b \subset\subset B \cap \Omega$, and let $E(B, b)$ be the connected component of $B \cap \Omega$ which contains $b$, i.e., $E(B, b) \supseteq b$. We are going to apply Theorem (i) of §2 with $\mathcal{V} = \bigcap_{p<q} \mathcal{O}L^p(\Omega)$ and $X = E(B, b)$. For this purpose we consider the linear operator

$$T : \bigcap_{p<q} \mathcal{O}L^p(\Omega) \to \mathbb{C}^{E(B,b)}, \ T(f)(z) := f(z) \ \text{for} \ z \in E(B, b).$$

For each fixed $z \in E(B, b)$, the fuctional

$$T_z : \bigcap_{p<q} \mathcal{OL}^p(\Omega) \to \mathbb{C}, \text{ defined by } T_z(f) = T(f)(z) = f(z), \text{ for } f \in \bigcap_{p<q} \mathcal{OL}^p(\Omega),$$

is continuous. (This follows from the fact that convergence in the space $\bigcap_{p<q} \mathcal{OL}^p(\Omega)$ implies pointwise convergence.) We also observe that, in this case, the set $S = \{f \in \mathcal{V} : T(f) \text{ is unbounded on } X\}$ is equal to

$$S(B,b) = \left\{ f \in \bigcap_{p<q} \mathcal{OL}^p(\Omega) : \sup_{z \in E(B,b)} |f(z)| = +\infty \right\}.$$

We claim that $S(B,b) \neq \varnothing$. Indeed, since the set $\overline{E(B,b)}$ meets the boundary of $\Omega$, there exists a point $\zeta \in \overline{E(B,b)} \cap \partial\Omega$. (See the remark in §3.) By hypotheses, there is a function $f_\zeta \in \bigcap_{p<q} \mathcal{OL}^p(\Omega)$ such that $\lim_{\substack{z \to \zeta \\ z \in \Omega}} f_\zeta(z) = \infty$, and, therefore $f_\zeta \in S(B,b)$. It follows (from Theorem (i) of §2) that $S(B,b)$ is dense and $\mathcal{G}_\delta$ in $\bigcap_{p<q} \mathcal{OL}^p(\Omega)$.

To complete the proof of the theorem, we consider a countable dense subset $\{w_1, w_2, w_3, \ldots\}$ of $\partial\Omega$, and the set $\mathcal{B} = \{B(w_j, \tau) : \tau \in \mathbb{Q}^+, \ j = 1,2,3,\ldots\}$. For each $B \in \mathcal{B}$, let $\Gamma_B$ be the countable set of the balls $b$ with centers in $(\mathbb{Q} + i\mathbb{Q})^n$ and rational radii, so that $b \subset\subset B \cap \Omega$. By Baire's theorem, the set

$$\bigcap_{B \in \mathcal{B}} \bigcap_{b \in \Gamma_B} S(B,b)$$

is dense and $\mathcal{G}_\delta$ in $\bigcap_{p<q} \mathcal{OL}^p(\Omega)$. Notice that if $f$ belongs to this set then $f$ is totally unbounded in $\Omega$. Indeed suppose that $w \in \partial\Omega$, $\varepsilon > 0$, and $E$ is a connected component of the set $B(w,\varepsilon) \cap \Omega$. Let $b$ be a ball with «rational» center and rational radius, and $b \subset\subset E$. Then we may choose a ball $B \in \mathcal{B}$ so that $B \subset B(w,\varepsilon)$ and $b \subset\subset B$. Then the connected component $E(B,b)$ of $B \cap \Omega$ which contains $b$, is contained in $E$, i.e., $E(B,b) \subset E$. Since $\sup_{z \in E(B,b)} |f(z)| = +\infty$, it follows that $\sup_{z \in E} |f(z)| = +\infty$.

To prove the last assertion of the theorem we will use Theorem (ii) of §2. For this purpose let us consider a pair of balls $(B,b)$ with $b \subset\subset B \cap \Omega \neq B$, and, as before, let $E(B,b)$ be the connected component of $B \cap \Omega$ which contains $b$. Then, by the remark in §3, $\overline{E(B,b)} \cap \partial\Omega \cap B \neq \varnothing$. If $\zeta \in \overline{E(B,b)} \cap \partial\Omega \cap B$, then the function $f_\zeta$ (of the hypotheses of the theorem) belongs to $\mathcal{V} = \bigcap_{p<q} \mathcal{OL}^p(\Omega)$ and its restriction $f_\zeta|_b$ (to $b$) does not have any bounded holomorphic extention to $B$. Hence Theorem (ii) of §2 gives the required conclusion.

## 5. Remarks.

**(i)** By examining the above proof we see that if the sets $B \cap \Omega$ are connected (for those $B$'s having sufficiently small radius – depending on the center of each $B$) then the theorem holds under the weaker hypothesis of the existence of the functions $f_\zeta$, not necessarily for all $\zeta \in \partial\Omega$, but only for $\zeta$ in a countable dense subset of $\partial\Omega$. This is the case – for example – in which the boundary of $\Omega$ is $C^1$.

**(ii)** Let us point out that the above theorem holds also in the case «$q = \infty$». The proof in this case is essentially the same. Although the case «$q = \infty$» is, in some sense, the most interesting one, it does not imply the case « $q < \infty$ ». Notices that changing the value of $q$ in $\bigcap\limits_{p<q} \mathcal{O}L^p(\Omega)$, changes not only the space but also the topology.

**(iii)** We can also prove an analogous theorem in the case of the spaces $\mathcal{O}L^p(\Omega)$ for each fixed $p$ ($1 \le p < \infty$). In this case we do not need to assume $\Omega$ to be bounded. Thus if $\Omega \subset \mathbb{C}^n$ is an open set and for every point $\zeta \in \partial\Omega$, there exists a function $f_\zeta$ such that $f_\zeta \in \mathcal{O}L^p(\Omega)$ and $\lim\limits_{z \in \Omega,\, z \to \zeta} f_\zeta(z) = \infty$, then the set of functions $g$ in the space $\mathcal{O}L^p(\Omega)$, which are totally unbounded in $\Omega$, is dense and $\mathcal{G}_\delta$ in this space.

**6. Theorem.** *Let $\Omega \subset \mathbb{C}^n$ be a bounded open set and $q \in \mathbb{R}$, $q > 1$. Suppose that for every point $\eta \in \partial\Omega$ and every $\varepsilon > 0$, there exists a function $f_{\eta,\varepsilon}$ such that*

$$f_{\eta,\varepsilon} \in \bigcap\limits_{p<q} \mathcal{O}L^p(\Omega) \ \ and \ \ f_{\eta,\varepsilon} \notin \mathcal{O}L^q(B(\zeta,\varepsilon) \cap \Omega).$$

*Then the set*

$$\mathcal{S}(\Omega, q) = \left\{ g \in \bigcap\limits_{p<q} \mathcal{O}L^p(\Omega) : g \notin \mathcal{O}L^q(B(\zeta,\varepsilon) \cap \Omega) \ for \ every \ \zeta \in \partial\Omega \ and \ every \ \varepsilon > 0 \right\}$$

*is dense and $\mathcal{G}_\delta$ in the space $\bigcap\limits_{p<q} \mathcal{O}L^p(\Omega)$.*

**Proof.** Let us fix a point $w \in \partial\Omega$ and $\varepsilon > 0$. We are going to apply Theorem (i) of §2 with $\mathcal{V} = \bigcap\limits_{p<q} \mathcal{O}L^p(\Omega)$ and $X$ being the set of all compact subsets $K$ of the intersection $B(w,\varepsilon) \cap \Omega$. For this purpose we consider the sublinear operator

$$T : \bigcap\limits_{p<q} \mathcal{O}L^p(\Omega) \to \mathbb{C}^X, \ T(f)(K) := \left( \int\limits_K |f|^q \, dv \right)^{1/q} \ \text{for } K \in X \, .$$

For every $K \in X$, the functional

$$T_K : \bigcap\limits_{p<q} \mathcal{O}L^p(\Omega) \to \mathbb{C}, \ T_K(f) = T(f)(K) \, ,$$

is continuous. Indeed, if $f_k$, $k = 1,2,3,\dots$, is a sequence which converges to $f$, in $\bigcap\limits_{p<q} \mathcal{O}L^p(\Omega)$, then $f_k$ converges to $f$, uniformly on $K$, and therefore

$$\int\limits_K |f_k|^q \, dv \to \int\limits_K |f|^q \, dv \text{, as } k \to \infty \, .$$

We also observe that, in this case, the set $S = \{f \in \mathcal{V} : T(f) \text{ is unbounded on } X\}$ is equal to

$$S(w,\varepsilon) = \left\{ f \in \bigcap\limits_{p<q} \mathcal{O}L^p(\Omega) : \int\limits_{B(w,\varepsilon) \cap \Omega} |f|^q \, dv = +\infty \right\} \, .$$

This follows from the fact that

$$\int\limits_{B(w,\varepsilon) \cap \Omega} |f|^q \, dv = \sup\limits_{K \in X} \int\limits_K |f|^q \, dv \, .$$

Also $S(w,\varepsilon) \ne \varnothing$, since $f_{w,\varepsilon} \in S(w,\varepsilon)$. Therefore, from Theorem (i) of §2, $S(w,\varepsilon)$ is dense and $\mathcal{G}_\delta$ in the space $\bigcap\limits_{p<q} \mathcal{O}L^p(\Omega)$.

Next let us observe that if $u_j$ is a sequence of points in $\partial\Omega$ which converges to a point $u \in \partial\Omega$, and $f \in \bigcap_{p<q} \mathcal{O}L^p(\Omega)$, then

$$\int_{B(u_j,\varepsilon)\cap\Omega} |f|^q \, dv = +\infty \ \ (\forall j) \ \Rightarrow \ \int_{B(u,2\varepsilon)\cap\Omega} |f|^q \, dv = +\infty \, .$$

This follows from the fact that if $\left| u_{j_0} - u \right| < \varepsilon$ (for some $j_0$) then $B(u,2\varepsilon) \supseteq B(u_{j_0},\varepsilon)$, which implies that

$$\int_{B(u,2\varepsilon)\cap\Omega} |f|^q \, dv \geq \int_{B(u_{j_0},\varepsilon)\cap\Omega} |f|^q \, dv \, .$$

To complete the proof of the theorem we consider a countable dense subset $\{w_1, w_2, w_3, \dots\}$ of $\partial\Omega$ and a decreasing sequence $\varepsilon_s$ of positive numbers, with $\varepsilon_s \to 0$. By the first part of the proof and Baire's theorem, the set

$$\bigcap_{j=1}^{\infty} \bigcap_{s=1}^{\infty} S(w_j, \varepsilon_s)$$

is dense and $\mathcal{G}_\delta$ in $\bigcap_{p<q} \mathcal{O}L^p(\Omega)$. Notice that if $f$ belongs to this set, and $\zeta \in \partial\Omega$, and $w_{j_m}$ is a subsequence of $w_j$ which converges to $\zeta$, then $\int_{B(w_{j_m},\varepsilon_s)\cap\Omega} |f|^q \, dv = +\infty$, and therefore $\int_{B(\zeta,2\varepsilon_s)\cap\Omega} |f|^q \, dv = +\infty$. Since this holds for every $\zeta \in \partial\Omega$, and the sequence $\varepsilon_s \to 0$, this implies that $f \in \mathcal{S}(\Omega, q)$. This completes the proof of the theorem.

**7. Remarks. (i)** By examining the above proof we see that this theorem holds under the weaker hypothesis of the existence of the functions $f_{\eta,\varepsilon}$, not necessarily for all $\eta \in \partial\Omega$, but only for $\eta$ in a countable dense subset of $\partial\Omega$.

**(ii)** The following version of the above theorem can be proved in a similar manner. Let $\Omega \subset \mathbb{C}^n$ be a bounded open set and $q, \tilde{q} \in \mathbb{R}$ with $\tilde{q} > q > 1$. Suppose that for every point $\zeta \in \partial\Omega$, there exists a function $f_\zeta$ such that

$$f_\zeta \in \bigcap_{p<q} \mathcal{O}L^p(\Omega) \ \text{ and } \ f_\zeta \notin \mathcal{O}L^{\tilde{q}}(B(\zeta,\varepsilon)\cap\Omega) \ \text{ for every } \ \varepsilon > 0 \, .$$

Then the set

$$\mathcal{S}(\Omega, q, \tilde{q}) = \left\{ g \in \bigcap_{p<q} \mathcal{O}L^p(\Omega) : g \notin \mathcal{O}L^{\tilde{q}}(B(\zeta,\varepsilon)\cap\Omega) \ \textit{for every } \zeta \in \partial\Omega \textit{ and every } \varepsilon > 0 \right\}$$

is dense and $\mathcal{G}_\delta$ in the space $\bigcap_{p<q} \mathcal{O}L^p(\Omega)$.

**(iii)** If the boundary of $\Omega$ is $C^1$ and a function $g \notin \mathcal{O}L^q(B(\zeta,\varepsilon)\cap\Omega)$, for every $\zeta \in \partial\Omega$ and every $\varepsilon > 0$, then $g$ is singular at every point of $\partial\Omega$. Indeed, this follows from the fact that for sufficiently small $\varepsilon > 0$ (depending on each point $\zeta \in \partial\Omega$), the sets $B(\zeta,\varepsilon)\cap\Omega$ are connected.

**(iv)** In the above theorem, if the sets $B \cap \Omega$ are connected (for those $B$'s having sufficiently small radius – depending on the center of each $B$) then the set of the functions $h$ in the space $\bigcap_{p<q} \mathcal{O}L^p(\Omega)$, which are singular at every boundary point of $\Omega$ is dense and $\mathcal{G}_\delta$ in this space. This follows from Theorem (ii) of §2.

**8. Theorem.** *Let* $\Omega \subset \mathbb{C}^n$ *be a bounded open set and* $1 < q \le \infty$. *Suppose that for every point* $\zeta \in \partial\Omega$, *there exists a function* $f_\zeta$ *such that*

$$f_\zeta \in \bigcap_{p<q} \mathcal{O}L^p(\Omega), \quad f_\zeta \notin \mathcal{O}L^q(B(\zeta,\varepsilon) \cap \Omega) \text{ for every } \varepsilon > 0, \text{ and } \lim_{\substack{z \to \zeta \\ z \in \Omega}} f_\zeta(z) = \infty.$$

*Then the set*

$$\left\{ g \in \bigcap_{p<q} \mathcal{O}L^p(\Omega) : g \text{ is totally unbounded in } \Omega \text{ and } g \notin \mathcal{O}L^q(B(\zeta,\varepsilon) \cap \Omega), \forall\, \zeta \in \partial\Omega \text{ and } \forall\, \varepsilon > 0 \right\}$$

*is dense and* $\mathcal{G}_\delta$ *in the space* $\bigcap_{p<q} \mathcal{O}L^p(\Omega)$.

**Proof.** The conclusion follows from Theorems 4 and 6. Indeed, it suffices to notice that the set in this theorem is the intersection of the corresponding sets of the Theorems 4 and 6, and that the intersection of two dense and $\mathcal{G}_\delta$ sets in the complete metric space $\bigcap_{p<q} \mathcal{O}L^p(\Omega)$ is again dense and $\mathcal{G}_\delta$, by Baire's theorem.

**9. Examples in the case** $n = 1$. **(i)** Let $\Omega \subset \mathbb{C}$ be a bounded open set with $C^1$ boundary. For a fixed point $\zeta \in \partial\Omega$, let us consider the holomorphic function

$$f_\zeta : \Omega \to \mathbb{C}, \ f_\zeta(z) = \frac{1}{z - \zeta}, \ z \in \Omega.$$

Then $f_\zeta \in \bigcap_{p<2} \mathcal{O}L^p(\Omega)$ but $f_\zeta \notin \mathcal{O}L^2(\Omega)$. Indeed, for «small» $\delta > 0$,

$$\iint_{z \in B(\zeta,\delta)} \frac{dv(z)}{|z-\zeta|^p} < +\infty \text{ when } p < 2, \text{ while } \iint_{z \in B(\zeta,\delta) \cap \Omega} \frac{dv(z)}{|z-\zeta|^2} = +\infty.$$

To prove the last equation, it suffices to notice that, since $\partial\Omega$ is assumed to be $C^1$, there is a small angular region A with vertex at $\zeta$ such that $A \cap B(\zeta,\delta) \subset \Omega$, and, that the integral

$$\iint_{z \in A \cap B(\zeta,\delta)} \frac{dv(z)}{|z-\zeta|^2} = +\infty,$$

as we can easily see if we integrate in polar coordinates with center at $\zeta$.

Next, continuing to consider the point $\zeta \in \partial\Omega$ fixed, let $a \in \mathbb{C} - \overline{\Omega}$ be a point, sufficiently close to the point $\zeta$, so that the line segment $[\zeta, a]$, which connects $a$ and $\zeta$, is contained in $\mathbb{C} - \Omega$. (Such a point exists since we assume that $\partial\Omega$ is $C^1$.) Then, in the set $\Omega$, there exists a holomorphic branch of $\log\left(\dfrac{z-a}{z-\zeta}\right)$, i.e., there exists a holomorphic function $g_\zeta(z)$, $z \in \Omega$, such that $\exp(g_\zeta(z)) = \dfrac{z-a}{z-\zeta}$. Indeed, the Möbius transformation $(z-a)/(z-\zeta)$ maps the point $a$ to $0$, $\zeta$ to $\infty$, and the line segment $[a,\zeta]$ to a half line in the complex plane, starting at $0$. We may also choose $g_\zeta$ so that $\left|\operatorname{Im} g_\zeta(z)\right| \le \pi$ for $z \in \Omega$. Then, for this function $g_\zeta$, the integral

$$(1) \qquad \iint_{z \in B(\zeta,\delta) \cap \Omega} \left| g_\zeta(z) \right|^p dv(z) < +\infty \text{ for every } p < \infty,$$

while $\lim_{z \in \Omega, \, z \to \zeta} g_\zeta(z) = \infty$. To prove (1), it suffices to notice that

$$(\log x)^p \le (k!)^{p/k} x^{p/k} \text{ for every } x > 1, \ p \ge 1 \text{ and } k \in \mathbb{N},$$

and that if $\log w = \log|w| + i\theta$, then

$$\left|\log w\right|^p = [(\log|w|)^2 + \theta^2]^{p/2}, \text{ for } \theta \in \mathbb{R}.$$

Indeed, since $g_\zeta(z) = \log\left|\dfrac{z-a}{z-\zeta}\right| + i\theta$ (with $|\theta| \le \pi$), it follows that, for $z \in \Omega$ which are sufficiently close to the point $\zeta$,

$$\left|g_\zeta(z)\right|^p = \left[\left(\log\left|\frac{z-a}{z-\zeta}\right|\right)^2 + \theta^2\right]^{p/2} \preceq (k!)^{p/k}\left|\frac{z-a}{z-\zeta}\right|^{p/k} \preceq (k!)^{p/k}\left|\frac{1}{z-\zeta}\right|^{p/k}.$$

Then (1) follows by an appropriate choice of $k \in \mathbb{N}$. Finally (1) implies that $g_\zeta \in \bigcap\limits_{p<\infty} \mathcal{O}L^p(\Omega)$, while $g_\zeta \notin \mathcal{O}L^\infty(\Omega)$.

**(ii)** With notation as in the previous example, and for $1 < q < \infty$, let us consider the function

$$h_{q,\zeta}(z) = \exp\left[\frac{2}{q} g_\zeta(z)\right], \ z \in \Omega.$$

Then $h_{q,\zeta} \in \bigcap\limits_{p<q} \mathcal{O}L^p(\Omega)$, while $h_{q,\zeta} \notin \mathcal{O}L^q(\Omega)$.

**(iii)** For $\alpha \in \mathbb{R}$, $\alpha \ge 1$, let us consider the domain $\Omega_\alpha = \{z = x + iy \in \mathbb{C} : 0 < x < 1 \ and \ 0 < y < x^\alpha\}$. Then

$$\frac{1}{z} \in \bigcap\limits_{p<\alpha+1} \mathcal{O}L^p(\Omega_\alpha) \text{ and } \frac{1}{z} \notin \mathcal{O}L^{\alpha+1}(\Omega_\alpha),$$

$$\log z \in \bigcap\limits_{p<\infty} \mathcal{O}L^p(\Omega_\alpha) \text{ and } \log z \notin \mathcal{O}L^\infty(\Omega_\alpha), \text{ and}$$

$$\frac{1}{z^{(\alpha+1)/q}} \in \bigcap\limits_{p<q} \mathcal{O}L^p(\Omega_\alpha) \text{ and } \frac{1}{z^{(\alpha+1)/q}} \notin \mathcal{O}L^q(\Omega_\alpha) \text{ for } q \in \mathbb{R}, \ q > 0.$$

**(iv)** Let $\Omega = \{z = x + iy \in \mathbb{C} : 0 < x < 1 \ and \ 0 < y < \exp(-1/x^2)\}$. Then

$$\frac{1}{z^N} \in \bigcap\limits_{p<\infty} \mathcal{O}L^p(\Omega) \text{ and } \frac{1}{z^N} \notin \mathcal{O}L^\infty(\Omega), \text{ for every } N \in \mathbb{N}.$$

**10. Theorem.** *(i) Let $\Omega \subset \mathbb{C}$ be an arbitrary bounded open set. Then the set of the functions $g \in \bigcap\limits_{p<2} \mathcal{O}L^p(\Omega)$ which are totally unbounded in $\Omega$ is dense and $\mathcal{G}_\delta$ in the space $\bigcap\limits_{p<2} \mathcal{O}L^p(\Omega)$.*

*(ii) Suppose that $\Omega \subset \mathbb{C}$ is a bouned open set such that for every point $\zeta \in \partial\Omega$, the connected component $C_\zeta$ of $\mathbb{C} - \Omega$ which contains $\zeta$, contains at least one more point, i.e., $C_\zeta - \{\zeta\} \neq \varnothing$. Then, for each fixed $q$ with $1 < q < \infty$, the set of the functions $g \in \bigcap\limits_{p<q} \mathcal{O}L^p(\Omega)$ which are totally unbounded in $\Omega$ is dense and $\mathcal{G}_\delta$ in the space $\bigcap\limits_{p<q} \mathcal{O}L^p(\Omega)$.*

*(iii) Suppose that $\Omega \subset \mathbb{C}$ is a bouned open set with $C^1$ boundary and $1 < q \le \infty$. Then the set*

$$\left\{g \in \bigcap\limits_{p<q} \mathcal{O}L^p(\Omega) : g \text{ is totally unbounded in } \Omega \text{ and } g \notin \mathcal{O}L^q(B(\zeta,\varepsilon) \cap \Omega), \forall \zeta \in \partial\Omega \text{ and } \forall \varepsilon > 0\right\}$$

*is dense and $\mathcal{G}_\delta$ in the space $\bigcap\limits_{p<q} \mathcal{O}L^p(\Omega)$.*

**Proof.** Having in mind the example (i) of §9, we easily obtain part (i), from Theorem 4, applied with the functions $\{f_\zeta : \zeta \in \partial\Omega\}$ where $f_\zeta : \Omega \to \mathbb{C}$, $f_\zeta(z) = 1/(z - \zeta)$, $z \in \Omega$. Part (iii) follows from Theorem 8, applied with the functions $g_\zeta$ of example (i) of §9 in the case $q = \infty$, and the functions $h_{q,\zeta}$ of example (ii) of §9 in the case $1 < q < \infty$. It remains to prove part (ii). For this purpose let us take a point $a \in C_\zeta$, $a \neq \zeta$, and define the set $K$ as follows: If $C_\zeta$ is unbounded then $K := C_\zeta \cup \{\infty\}$, and if $C_\zeta$ is bounded then $K := C_\zeta$. In both cases $K$ is a compact and connected set (in $\mathbb{C} \cup \{\infty\}$) containing the points $\zeta$ and $a$. Then the Möbius transformation $(z - a)/(z - \zeta)$ maps the point $a$ to $0$, $\zeta$ to $\infty$, and the set $K$ to a connected compact set $\Gamma$ containing the points $0$ and $\infty$. Then, in the open simply connected set $\mathbb{C} - \Gamma$, there is a holomorphic branch of the logarithm, and, therefore, there is a function $\varphi_\zeta(z)$, holomorphic in $z \in \Omega$, such that $\exp[\varphi_\zeta(z)] = (z - a)/(z - \zeta)$. Also the function

$$\psi_\zeta(z) := \exp\left[\frac{2}{q} \varphi_\zeta(z)\right]$$

is holomorphic in $\Omega$ and

$$\left|\psi_\zeta(z)\right|^p = \exp\left[\frac{2p}{q} \operatorname{Re} \varphi_\zeta(z)\right] = \left|\frac{z - a}{z - \zeta}\right|^{2p/q}.$$

Therefore $\psi_\zeta \in \bigcap_{p < q} \mathcal{O}L^p(\Omega)$, and, since $\lim_{z \to \zeta,\, z \in \Omega} \psi_\zeta(z) = \infty$, part (ii) follows from Theorem 4.

## 11. The case of the unit ball of $\mathbb{C}^n$.

Let us consider the unit ball $\Omega = \{z \in \mathbb{C}^n : |z| < 1\}$. For fixed $\zeta \in \partial\Omega$, we consider the function

$$f_\zeta(z) = \frac{1}{1 - \langle z, \zeta \rangle} = \frac{1}{1 - \sum_{j=1}^n \overline{\zeta}_j z_j}, \quad z \in \Omega.$$

Then

$$f_\zeta \in \bigcap_{p < n+1} \mathcal{O}L^p(\Omega) \text{ and } f_\zeta \notin \mathcal{O}L^{n+1}(\Omega).$$

Indeed, if $p < n + 1$ then the integral

$$\int_\Omega \frac{dv(z)}{\left|1 - \langle z, \eta \rangle\right|^p},$$

as a function of $\eta$, remains bounded for $\eta \in \Omega$ (see [14] and [16]), and, therefore, letting $\eta \to \zeta$,

$$\int_\Omega \frac{dv(z)}{\left|1 - \langle z, \zeta \rangle\right|^p} = \int_\Omega \lim_{\eta \to \zeta} \frac{dv(z)}{\left|1 - \langle z, \eta \rangle\right|^p} \leq \liminf_{\eta \to \zeta} \int_\Omega \frac{dv(z)}{\left|1 - \langle z, \eta \rangle\right|^p} < +\infty.$$

Next we show that

(1)
$$\int_\Omega \frac{dv(z)}{\left|1 - \langle z, \zeta \rangle\right|^{n+1}} = +\infty.$$

Indeed, for $r < 1$ (sufficiently close to 1),

$$\int_\Omega \frac{dv(z)}{\left|1 - \langle z, r\zeta \rangle\right|^{n+1}} \geq \lambda \log \frac{1}{1 - r^2},$$

where $\lambda$ is a positive constant independant of $r$. (See [14] and [16].) Since

$$\int_\Omega \frac{dv(z)}{\left|1 - \langle z, r\zeta \rangle\right|^{n+1}} = \int_\Omega \frac{dv(z)}{\left|1 - \langle rz, \zeta \rangle\right|^{n+1}} = \frac{1}{r^{2n}} \int_{r\Omega} \frac{dv(z)}{\left|1 - \langle z, \zeta \rangle\right|^{n+1}} \quad (\text{where } r\Omega = \{z \in \mathbb{C}^n : |z| < r\}),$$

it follows that

$$\int_{r\Omega} \frac{dv(z)}{\left|1 - \langle z, \zeta \rangle\right|^{n+1}} \geq \lambda r^{2n} \log \frac{1}{1 - r^2}.$$

Letting $r \to 1^-$, we obtain (1).

Observing that $\mathrm{Re}(1 - \langle z, \zeta \rangle) > 0$, for $z \in \Omega$, we see that $\mathrm{Re}\, f_\zeta(z) > 0$, and therefore $\log f_\zeta(z)$ is defined and holomorphic for $z \in \Omega$, where $\log$ is the principal branch of the logarithm with $\left|\arg\right| \leq \pi$. Also $\left|\mathrm{Im}[\log f_\zeta(z)]\right| \leq \pi/2$. It follows, as in example (i) of §9, that

$$\log f_\zeta \in \bigcap_{p < \infty} \mathcal{O}L^p(\Omega), \text{ while } \log f_\zeta \notin \mathcal{O}L^\infty(\Omega).$$

Also the function $\left(f_\zeta\right)^{(n+1)/q} = \exp\left[\frac{n+1}{q} \log f_\zeta\right]$ satisfies

$$\left(f_\zeta\right)^{(n+1)/q} \in \bigcap_{p < q} \mathcal{O}L^p(\Omega) \text{ and } \left(f_\zeta\right)^{(n+1)/q} \notin \mathcal{O}L^q(\Omega) \text{ for } q \in \mathbb{R}, \ q > 0.$$

**12. Theorem.** *Let* $1 < q \leq \infty$. *If* $\Omega$ *is the unit ball of* $\mathbb{C}^n$, *or more generally if* $\Omega$ *is an ellipsoid as in* §13(ii), *below, then the set*

$$\left\{ g \in \bigcap_{p < q} \mathcal{O}L^p(\Omega) : g \text{ is totally unbounded in } \Omega \text{ and } g \notin \mathcal{O}L^q(B(\zeta, \varepsilon) \cap \Omega), \forall \zeta \in \partial\Omega \text{ and } \forall \varepsilon > 0 \right\}$$

*is dense and* $\mathcal{G}_\delta$ *in the space* $\bigcap\limits_{p < q} \mathcal{O}L^p(\Omega)$.

**Proof.** It suffices to apply Theorem 8, with appropriate choices from the set of the functions which were constructed in §11 and §13(ii).

**13. The case of convex sets. (i)** Let $\Omega \subset \mathbb{C}^n$ be a bounded open and convex set with $C^1$ boundary, and let $\rho$ be a defining function for $\Omega$. For a fixed point $\zeta \in \partial\Omega$, we consider the function

$$(1) \qquad f_\zeta(z) = \frac{1}{\sum\limits_{j=1}^{n} \dfrac{\partial \rho}{\partial \zeta_j}(\zeta)(\zeta_j - z_j)}, \ z \in \Omega.$$

Then

$$(2) \qquad f_\zeta \in \bigcap_{p < 2} \mathcal{O}L^p(\Omega), \text{ while } f_\zeta \notin \mathcal{O}L^\infty(\Omega), \text{ and}$$

$$(3) \qquad \log f_\zeta \in \bigcap_{p < \infty} \mathcal{O}L^p(\Omega), \text{ while } \log f_\zeta \notin \mathcal{O}L^\infty(\Omega).$$

To prove (2), we will show that for $p < 2$,

$$(4) \qquad \int_{B(\zeta, \delta) \cap \Omega} \left|f_\zeta(z)\right|^p dv(z) < +\infty, \text{ for «small» } \delta > 0.$$

Assuming, without loss of generality, that $\dfrac{\partial \rho}{\partial z_1}(\zeta) \neq 0$, let us consider the $\mathbb{C}$-affine transformation

$$w_1(z) = \sum_{j=1}^{n} \frac{\partial \rho}{\partial \zeta_j}(\zeta)(\zeta_j - z_j), \ w_2(z) = \zeta_2 - z_2, \ \ldots, \ w_n(z) = \zeta_n - z_n.$$

Using this transformation we see that (2) follows from the fact that

$$\int_{|w| < \tilde{\delta}} \frac{dv(w)}{\left|w_1\right|^p} < +\infty \ \text{(for } \tilde{\delta} > 0\text{)}.$$

To justify (3), let us recall that the convexity of $\Omega$ implies that

$$\text{Re}\left[\sum_{j=1}^{n}\frac{\partial\rho}{\partial\zeta_j}(\zeta)(\zeta_j-z_j)\right]>0 \ \text{ for every } \ z\in\Omega\,,$$

so that $\log f_\zeta$ is well defined and holomorphic in $\Omega$. Then, using (3) as in example (i) of §9, we see that, for «small» $\delta>0$,

$$\int\limits_{B(\zeta,\delta)\cap\Omega}\left|\log f_\zeta(z)\right|^p dv(z)<+\infty\,, \text{ for every } \ p<\infty\,,$$

and this implies (3).

We point out that in general the conclusion $f_\zeta\in\bigcap_{p<2}\mathcal{O}L^p(\Omega)$ cannot be improved in the sense that in some cases

$$\int\limits_{B(\zeta,\delta)\cap\Omega}\left|f_\zeta(z)\right|^2 dv(z)=+\infty$$

(see the example (iii) below).

We also observe that if $\lambda$ is another defining function for $\Omega$, then $\lambda=\varphi\rho$ for a continuous and positive function $\varphi$, and

$$\sum_{j=1}^{n}\frac{\partial\lambda}{\partial\zeta_j}(\zeta)(\zeta_j-z_j)=\varphi(\zeta)\sum_{j=1}^{n}\frac{\partial\rho}{\partial\zeta_j}(\zeta)(\zeta_j-z_j)\,.$$

Thus the function (1) does not depend «essentially» on the choice of the defining function of $\Omega$, as long as this set is convex with $C^1$ boundary. Notice also that the functions $f_\zeta$ depend continuously on the point $\zeta$.

**(ii)** Let $\Omega\subset\mathbb{C}^n$ be an ellipsoid of the form

$$\Omega=\left\{z\in\mathbb{C}^n:\frac{|z_1|^2}{a_1^2}+\frac{|z_2|^2}{a_2^2}+\cdots+\frac{|z_n|^2}{a_n^2}<1\right\},$$

where $a_1,a_2,...,a_n$ are some positive numbers. Using the defining function

$$\rho(z)=\frac{z_1\bar z_1}{a_1^2}+\frac{z_2\bar z_2}{a_2^2}+\cdots+\frac{z_n\bar z_n}{a_n^2}-1\,,$$

we find, for a fixed point $\zeta\in\partial\Omega$, that the function

$$f_\zeta(z)=\frac{1}{\sum_{j=1}^{n}\dfrac{\partial\rho}{\partial\zeta_j}(\zeta)(\zeta_j-z_j)}=\frac{1}{1-\sum_{j=1}^{n}\dfrac{\bar\zeta_j z_j}{a_j^2}} \ \text{ is defined and holomorphic for } z\in\Omega\,.$$

Then

$$\int\limits_{\Omega}\left|f_\zeta(z)\right|^p dv(z)<+\infty \ \text{ for every } \ p<n+1\,, \ \text{ and } \int\limits_{\Omega}\left|f_\zeta(z)\right|^{n+1} dv(z)=+\infty\,.$$

The above assertions can easily be reduced to the case of the unit ball, which we studied in §11. In fact it suffices to change the variables by setting $w_j=z_j/a_j$ and $\eta_j=\zeta_j/a_j$.

**(iii)** Let us consider the convex domain $D=\{z=(z_1,...,z_n)\in\mathbb{C}^n:|z|<1 \ and \ \text{Re } z_1>0\}$ and, as local defining function for $D$ near its boundary point $\zeta=0$ $(0\in\partial D)$, $\rho(z)=-(z_1+\bar z_1)/2$. Then the function (1) becomes $f_\zeta(z)=1/z_1$. In this case

$$\int\limits_{B(\zeta,\delta)\cap D}\left|f_\zeta(z)\right|^2 dv(z)=\int\limits_{B(0,\delta)\cap D}\frac{1}{|z_1|^2}\,dv(z)=+\infty\,, \text{ for every } \ \delta>0\,.$$

A similar computation can be done for every point $\zeta$ in the part of the boundary of $\partial D$ where $\operatorname{Re}\zeta = 0$ (and $|\zeta| < 1$). Of course at the points $\zeta \in \partial D$ where $\operatorname{Re}\zeta > 0$, the corresponding function $f_\zeta$ satisfies $f_\zeta \in \bigcap\limits_{p < n+1} \mathcal{O}L^p(D)$ and $f_\zeta \notin \mathcal{O}L^{n+1}(D)$, as we proved in §11.

**(iv)** Similarly to the previous example, if
$$R = \{z = (z_1,...,z_n) \in \mathbb{C}^n : 0 < \operatorname{Re} z_j < 1 \ and \ 0 < \operatorname{Im} z_j < 1, \ j = 1,2,...,n\},$$
then for every point $\zeta \in \partial R$ (where $\partial R$ is smooth), the function $f_\zeta$ satisfies
$$f_\zeta \in \bigcap_{p<2} \mathcal{O}L^p(R) \ \text{and} \ f_\zeta \notin \mathcal{O}L^2(R).$$
Similar conclusions hold for «most» points in the boundary of the polydisk
$$P = \{z = (z_1,...,z_n) \in \mathbb{C}^n : |z_j| < 1, \ j = 1,2,...,n\}.$$

**(v)** The results of example (i) can easily be extended to arbitrary convex sets – with no smoothness of its boundary assumed. More precisely let us consider any bounded open and convex set $\Omega \subset \subset \mathbb{C}^n$, and let us fix a point $\zeta \in \partial\Omega$. By the convexity of $\Omega$, there exist real numbers $\alpha_j = \alpha_j(\zeta)$, $\beta_j = \beta_j(\zeta)$, $j = 1,2,...,n$, such that $\sum\left[\left|\alpha_j\right|^2 + \left|\beta_j\right|^2\right] \neq 0$ and

$$\sum_{j=1}^{n}\left\{\alpha_j[x_j(z) - x_j(\zeta)] + \beta_j[y_j(z) - y_j(\zeta)]\right\} > 0 \ \text{for every} \ z \in \Omega,$$

where $x_j(z) = \operatorname{Re} z_j$, $y_j(z) = \operatorname{Im} z_j$, $x_j(\zeta) = \operatorname{Re}\zeta_j$, $y_j(\zeta) = \operatorname{Im}\zeta_j$. Setting $c_j := \alpha_j - i\beta_j$, we obtain
$$\operatorname{Re}\left[\sum_{j=1}^{n} c_j(z_j - \zeta_j)\right] > 0 \ \text{for every} \ z \in \Omega.$$
Then the conclusions of example (i) hold for the function $f_\zeta$ where
$$f_\zeta(z) = \frac{1}{\sum\limits_{j=1}^{n} c_j(z_j - \zeta_j)}, \ z \in \Omega.$$
For example the function $\log f_\zeta(z)$, $z \in \Omega$, belongs to the space $\bigcap\limits_{p<\infty} \mathcal{O}L^p(\Omega)$ and $\lim\limits_{\substack{z \to \zeta \\ z \in \Omega}} f_\zeta(z) = \infty$.

**14. Theorem.** *Let $\Omega \subset \subset \mathbb{C}^n$ be any bounded open and convex set and $1 < q \leq \infty$. Then the set of the functions $g$ in $\bigcap\limits_{p<q} \mathcal{O}L^p(\Omega)$ such that $g$ is totally unbounded in $\Omega$, is dense and $\mathcal{G}_\delta$ in the space $\bigcap\limits_{p<q} \mathcal{O}L^p(\Omega)$.*

**Proof.** It follows from Theorem 4 applied with the functions $\log f_\zeta$ of the above example (v).

**15. The case of strictly pseudoconvex domains.** In this section we will show that some functions which are defined in terms of Henkin's support function belong to certain Bergman spaces. First we describe Henkin's support function $\Phi(z,\zeta)$ which is constructed in [5]. Following Henkin and Leiterer [5], let us consider an open set $\Theta \subset \subset \mathbb{C}^n$ and a $C^2$ strictly plurisubharmonic function $\rho$ in a neighbourhood of $\overline{\Theta}$. If we set

$$\beta = \frac{1}{3}\min\left\{\sum_{1 \leq j,k \leq n} \frac{\partial^2 \rho(\zeta)}{\partial \zeta_j \partial \overline{\zeta}_k}\xi_j\overline{\xi}_k : \zeta \in \overline{\Theta}, \ \xi \in \mathbb{C}^n \ with \ |\xi| = 1\right\}$$

then $\beta > 0$ and there exist $C^1$ functions $a_{jk}$ in a neighbourhood of $\overline{\Theta}$ such that

$$\max\left\{\left|a_{jk}(\zeta) - \frac{\partial^2 \rho(\zeta)}{\partial \zeta_j \partial \zeta_k}\right| : \zeta \in \overline{\Theta}\right\} < \frac{\beta}{n^2}.$$

Let $\varepsilon > 0$ be sufficiently small so that

$$\max\left\{\left|\frac{\partial^2 \rho(\zeta)}{\partial x_j \partial x_k} - \frac{\partial^2 \rho(z)}{\partial x_j \partial x_k}\right| : \zeta, z \in \overline{\Theta} \ \text{with} \ |\zeta - z| \le \varepsilon\right\} < \frac{\beta}{2n^2} \ \text{for} \ j, k = 1, 2, \ldots, 2n,$$

where $x_j = x_j(\xi)$ are the real coordinates of $\xi \in \mathbb{C}^n$ such that $\xi_j = x_j(\xi) + ix_{j+n}(\xi)$. For $z, \zeta \in \overline{\Theta}$ we consider the modified Levi polynomial

$$Q(z, \zeta) = -\left[2\sum_{j=1}^n \frac{\partial \rho(\zeta)}{\partial \zeta_j}(z_j - \zeta_j) + \sum_{1 \le j, k \le n} a_{jk}(z_j - \zeta_j)(z_k - \zeta_k)\right].$$

Then we have the estimate

$$\text{Re} \, Q(z, \zeta) \ge \rho(\zeta) - \rho(z) + \beta|\zeta - z|^2 \ \text{for} \ z, \zeta \in \overline{\Theta} \ \text{with} \ |\zeta - z| \le \varepsilon.$$

The following theorem is proved in [5].

**Theorem.** *Let $\Omega \subset\subset \mathbb{C}^n$ be a strictly pseudoconvex open set, let $\Theta$ be an open neighbourhood of $\partial\Omega$, and let $\rho$ be a $C^2$ strictly plurisubharmonic function in a neighbourhood of $\overline{\Theta}$ such that*
$$\Omega \cap \Theta = \{z \in \Theta : \rho(z) < 0\}.$$
*Let us choose $\varepsilon$, $\beta$, and $P(\zeta, \zeta)$, as above, and let us make the positive number $\varepsilon$ smaller so that*
$$\{z \in \mathbb{C}^n : |\zeta - z| \le 2\varepsilon\} \subseteq \Theta \ \text{for every} \ \zeta \in \partial\Omega.$$
*Then there exists a function $\Phi(z, \zeta)$ defined for $\zeta$ in some open neighbourhood $U_{\partial\Omega} \subseteq \Theta$ of $\partial\Omega$ and $z \in U_{\overline{\Omega}} = \Omega \cup U_{\partial\Omega}$, which is $C^1$ in $(z, \zeta) \in U_{\overline{\Omega}} \times U_{\partial\Omega}$, holomorphic in $z \in U_{\overline{\Omega}}$, and such that $\Phi(z, \zeta) \ne 0$ for $(z, \zeta) \in U_{\overline{\Omega}} \times U_{\partial\Omega}$ with $|z - \zeta| \ge \varepsilon$, and*
$$\Phi(z, \zeta) = Q(z, \zeta) C(z, \zeta) \ \text{for} \ (z, \zeta) \in U_{\overline{\Omega}} \times U_{\partial\Omega} \ \text{with} \ |z - \zeta| \le \varepsilon,$$
*for some $C^1$−function $C(z, \zeta)$ defined for $(z, \zeta) \in U_{\overline{\Omega}} \times U_{\partial\Omega}$ and $\ne 0$ when $|z - \zeta| \le \varepsilon$.*

In this setting we will prove the following proposition. We use a set of coordinates – *the Levi coordinates* – which are appropriate when we are dealing with integrals involving the function $\Phi(z, \zeta)$. (See [5] and [13].) As a matter of fact we will use a slight modification of the Levi coordinates.

**Proposition.** *If, in addition, $\partial\Omega$ is $C^1$, then, for each fixed $\zeta \in \partial\Omega$ and for every $\delta > 0$,*
$$\int_{z \in B(\zeta, \delta) \cap \Omega} \frac{dv(z)}{|\Phi(z, \zeta)|^p} < +\infty \ \text{when} \ p < n+1, \ \text{and} \ \int_{z \in B(\zeta, \delta) \cap \Omega} \frac{dv(z)}{|\Phi(z, \zeta)|^{2n}} = +\infty.$$
*Therefore $\dfrac{1}{\Phi(\cdot, \zeta)} \in \bigcap_{p < n+1} \mathcal{O}L^p(\Omega)$ and $\dfrac{1}{\Phi(\cdot, \zeta)} \notin \mathcal{O}L^{2n}(\Omega)$. Furthermore the functions $\dfrac{1}{\Phi(z, \zeta)}$ are $C^1$ in $\zeta$.*

**Proof.** Since we assume $\partial\Omega$ to be $C^1$, $\nabla\rho \ne 0$ at the points of $\partial\Omega$. Having fixed $\zeta \in \partial\Omega$, we consider a coordinate system $t = (t_1, t_2, t_3, \ldots, t_{2n}) = (t_1(z), t_2(z), t_3(z), \ldots, t_{2n}(z))$, of real $C^1$−functions, for points $z \in \mathbb{C}^n = \mathbb{R}^{2n}$, which are sufficiently close to the point $\zeta$, as follows: We set
$$t_1(z) = -\rho(z) \ \text{and} \ t_2(z) = \text{Im} \, Q(z, \zeta).$$

Then $d_z Q(z,\zeta)\big|_{z=\zeta} = -2\sum_{j=1}^{n} \dfrac{\partial \rho(\zeta)}{\partial \zeta_j}\,dz_j \bigg|_{z=\zeta} = -2\partial\rho(\zeta)$ and, therefore,

$$d_z t_2(z)\big|_{z=\zeta} = d_z[\operatorname{Im} Q(z,\zeta)]\big|_{z=\zeta} = i[\partial\rho(\zeta) - \bar\partial\rho(\zeta)].$$

On the other hand,

$$d_z t_1(z)\big|_{z=\zeta} = d_z[-\rho(z)]\big|_{z=\zeta} = -[\partial\rho(\zeta) + \bar\partial\rho(\zeta)].$$

It follows that

$$\big(d_z t_1(z)\big|_{z=\zeta}\big) \wedge \big(d_z t_2(z)\big|_{z=\zeta}\big) = -2i\,\partial\rho(\zeta) \wedge \bar\partial\rho(\zeta) \neq 0.$$

Now the existence of $C^1-$functions $t_3(z),...,t_{2n}(z)$ such that the mapping

$$z \to (t_1(z), t_2(z), t_3(z),..., t_{2n}(z))$$

is a $C^1-$diffeomorphism, from an open neighbourhood of the point $\zeta$ to an open neighbourhood of $0 \in \mathbb{C}^n = \mathbb{R}^{2n}$ (with $t(\zeta)=0$), follows from the inverse function theorem. Also let us point out that, for $z$ sufficiently close to $\zeta$, $z \in \Omega$ if and only if $t_1 = -\rho(z) > 0$.

We will show that, for every $\delta > 0$,

(1)     $\displaystyle\int_{z \in B(\zeta,\delta)\cap\Omega} \dfrac{dv(z)}{|\Phi(z,\zeta)|^p} < +\infty \quad$ for $\ p < n+1$.

For points $z \in \Omega$ which are sufficiently close to $\zeta$,

$$|\Phi(z,\zeta)| \approx |Q(z,\zeta)| \approx |\operatorname{Re} Q(z,\zeta)| + |\operatorname{Im} Q(z,\zeta)| \geq -\rho(z) + \beta|\zeta - z|^2 + |\operatorname{Im} Q(z,\zeta)|$$

and

$$|\zeta - z|^2 \approx t_1^2 + t_2^2 + t_3^2 + \cdots + t_{2n}^2.$$

(When we write $\mathrm{A} \approx \mathrm{B}$, we mean that $\lambda\mathrm{B} \leq \mathrm{A} \leq \mu\mathrm{B}$, for some positive constants $\lambda$ and $\mu$ which are independent of $z$.)

Therefore (for $z \in \Omega$ and sufficiently close to $\zeta$)

$$|\Phi(z,\zeta)| \succeq t_1 + t_1^2 + t_2^2 + t_3^2 + \cdots + t_{2n}^2 + |t_2|.$$

(When we write $\mathrm{A} \succeq \mathrm{B}$, we mean that $\mathrm{A} \geq \lambda\mathrm{B}$, for some positive constant $\lambda$ which is independent of $z$.) Therefore (1) follows from

$$\int_{t_1 > 0} \dfrac{dt}{(t_1 + |t_2| + t_1^2 + t_2^2 + t_3^2 + \cdots + t_{2n}^2)^p} < +\infty$$

or equivalently from

$$\int_{t_1 > 0} \dfrac{dt}{(t_1 + |t_2| + t_3^2 + \cdots + t_{2n}^2)^p} < +\infty \quad (\ p < n+1).$$

(In the above integrals $dt = dt_1 dt_2 \cdots dt_{2n}$ and $t$ is restricted in a «small» neighbourhood of $0 \in \mathbb{C}^n = \mathbb{R}^{2n}$, i.e., $|t|$ is «small».)

We will also show that, for every $\delta > 0$,

(2)     $\displaystyle\int_{z \in B(\zeta,\delta)\cap\Omega} \dfrac{dv(z)}{|\Phi(z,\zeta)|^{2n}} = +\infty$.

This time we will use the fact that, for points $z \in \Omega$ which are sufficiently close to $\zeta$,

$$|\Phi(z,\zeta)| \approx |Q(z,\zeta)| \preceq |\zeta - z| \approx (t_1^2 + t_2^2 + t_3^2 + \cdots + t_{2n}^2)^{1/2}.$$

Therefore (2) follows from

$$\int_{t_1 > 0} \dfrac{dv(z)}{(t_1^2 + t_2^2 + t_3^2 + \cdots + t_{2n}^2)^n} = +\infty.$$

This completes the proof of the proposition.

**16. Theorem.** *Let $\Omega \subset \mathbb{C}^n$ be a strictly pseudoconvex open set with $C^2$ − boundary, and $1 < q < \infty$. Then the following hold:*

*(i) For every point $\zeta \in \partial\Omega$, there exists a function $f_\zeta$ such that*

$$f_\zeta \in \bigcap_{p<\infty} \mathcal{O}L^p(\Omega) \quad and \quad \lim_{\substack{z \to \zeta \\ z \in \Omega}} f_\zeta(z) = \infty.$$

*(ii) For every point $\zeta \in \partial\Omega$, there exists a function $h_\zeta$ such that*

$$h_\zeta \in \bigcap_{p<q} \mathcal{O}L^p(\Omega) \ and \ h_\zeta \notin \mathcal{O}L^{2nq/(n+1)}(B(\zeta,\delta)\cap\Omega) \ for \ every \ \delta>0, \ and \ \lim_{\substack{z \to \zeta \\ z \in \Omega}} h_\zeta(z)=\infty.$$

*(iii) The set*

$$\left\{ g \in \bigcap_{p<\infty} \mathcal{O}L^p(\Omega) : g \ is \ totally \ unbounded \ in \ \Omega \right\}$$

*is dense and $\mathcal{G}_\delta$ in the space $\bigcap\limits_{p<\infty} \mathcal{O}L^p(\Omega)$.*

*(iv) The set*

$$\left\{ g \in \bigcap_{p<q} \mathcal{O}L^p(\Omega) : g \ is \ totally \ unbounded \ in \ \Omega \ and \ g \notin \mathcal{O}L^{2nq/(n+1)}(B(\zeta,\delta)\cap\Omega) , \forall \zeta \in \partial\Omega \ and \ \forall \delta > 0 \right\}$$

*is dense and $\mathcal{G}_\delta$ in the space $\bigcap\limits_{p<q} \mathcal{O}L^p(\Omega)$.*

**Proof.** Let $\rho$ be a $C^2$ strictly plurisubharmonic defining function of $\Omega$, defined in an open neighbourhood of $\overline{\Omega}$. Let us also fix a point $\zeta \in \partial\Omega$. Then, as it follows from Taylor's theorem and the strict plurisubharmonicity of $\rho$ (see [13]), the Levi polynomial of $\rho$

$$F(z,\zeta) = -\left[ 2\sum_{j=1}^n \frac{\partial\rho(\zeta)}{\partial\zeta_j}(z_j - \zeta_j) + \sum_{1\le j,k\le n} \frac{\partial^2\rho(\zeta)}{\partial\zeta_j\partial\zeta_k}(z_j-\zeta_j)(z_k-\zeta_k) \right]$$

satisfies the inequality

$$\operatorname{Re} F(z,\zeta) \ge \rho(\zeta) - \rho(z) + \beta|\zeta - z|^2 \ \text{for} \ z \in \mathbb{C}^n \ \text{with} \ |\zeta - z| < \varepsilon,$$

for some «small» positive constants $\varepsilon$ and $\beta$. In particular,

$$\operatorname{Re} F(z,\zeta) > 0 \ \text{for} \ z \in B(\zeta,\varepsilon) \cap \overline{\Omega} - \{\zeta\}.$$

It follows that the function $\log[1/F(z,\zeta)]$ is defined and holomorphic for $z \in B(\zeta,\varepsilon) \cap \Omega$, and that $\lim\limits_{z\in\Omega, z\to\zeta} \log[1/F(z,\zeta)] = \infty$. (Here $\log$ is the principal branch of the logarithm with $|\arg| \le \pi$.) Also we can prove, as in the proof of the proposition of §15, that if $q < n+1$,

$$(1) \qquad \int_{z\in B(\zeta,\delta)\cap\Omega} \frac{dv(z)}{|F(z,\zeta)|^q} < +\infty \ \text{for every} \ \delta > 0.$$

Then, using (1) (with $q = 1$, for example) as in example (i) of §9, we obtain

$$(2) \qquad \int_{B(\zeta,2\varepsilon/3)\cap\Omega} \left| \log\left[\frac{1}{F(z,\zeta)}\right] \right|^p dv(z) < +\infty, \ \text{for every} \ p < \infty.$$

Next we consider a $C^\infty$ − function $\chi : \mathbb{C}^n \to \mathbb{R}$, $0 \le \chi(z) \le 1$, with compact support contained in $B(\zeta, 2\varepsilon/3)$, and such that $\chi(z) = 1$ when $z \in B(\zeta, \varepsilon/3)$. Now the function

$$\chi(z) \log\left[\frac{1}{F(z,\zeta)}\right]$$

is extended to a $C^\infty$ − function in $\Omega$, by defining it to be $0$ in $\Omega - B(\zeta, 2\varepsilon/3)$. Then the $(0,1)$ − form

$$u(z) := \overline{\partial} \left\{ \chi(z) \log \left[ \frac{1}{F(z,\zeta)} \right] \right\}$$

is defined and is $C^\infty$ in a open neighbourhood $\overline{\Omega}$, it is zero for $z \in B(\zeta, \varepsilon/3) \cap \Omega$, and, in particular, it has bounded coefficients in $\Omega$. In fact $u(z)$ extends to a $C^\infty$ $(0,1)$ − form for $z$ in an open neighbourhood of $\overline{\Omega}$, since the function $\log \left[ \dfrac{1}{F(z,\zeta)} \right]$ is holomorphic in an open neighbourhood of the compact set $\overline{[B(\zeta, 2\varepsilon/3) - B(\zeta, \varepsilon/3)] \cap \Omega}$. It follows that there exists a bounded $C^\infty$ − function $\psi : \Omega \to \mathbb{C}$ which solves the equation $\overline{\partial} \psi = u$ in $\Omega$. (See [13].) Then the function

$$f_\zeta(z) := \chi(z) \log \left[ \frac{1}{F(z,\zeta)} \right] - \psi(z)$$

satisfies the requirements of (i) (as it follows from (2)).

A function $h_\zeta$ which satisfies the requirements of (ii) is

$$h_\zeta(z) = \exp \left[ \frac{n+1}{q} f_\zeta(z) \right] = \exp \left\{ \chi(z) \log \left[ \frac{1}{F(z,\zeta)} \right]^{(n+1)/q} - \frac{n+1}{q} \psi(z) \right\}.$$

Indeed, we have

$$\int_{z \in B(\zeta,\delta) \cap \Omega} \frac{dv(z)}{|F(z,\zeta)|^{2n}} = +\infty \quad \text{(for every } \delta > 0 \text{)}$$

(this is proved in the same manner as the analogous result of the Proposition in §15) which implies that

$$\int_{z \in B(\zeta,\delta) \cap \Omega} \frac{dv(z)}{|h_\zeta(z)|^{2nq/(n+1)}} = +\infty.$$

Notice that the behaviour of the above integral is not affected by the functions $\chi$ or $\psi$, since $\chi \equiv 1$ near $\zeta$ and $\psi$ is bounded in $\Omega$ (so that $\exp[-(n+1)\psi/q]$ is both bounded and bounded away from zero in $\Omega$).

Finally assertions (iii) and (iv) follow from (i) and (ii), in combination with Theorems 4 and 7. (See also the remark (ii) of §7.)

**17. Remark.** It follows from the above theorem, in combination with Theorem (ii) of §2, that the set of the functions $h$ in the space $\bigcap_{p<q} \mathcal{O}L^p(\Omega)$ which are singular at every boundary point of $\partial\Omega$ is dense and $\mathcal{G}_\delta$ in this space, for $1 < q \leq \infty$. (Similar conclusions are reached also in the case of the convex domains, following Theorem 14, and certain − more general − domains in $\mathbb{C}$, following Theorem 10.)

**18. The case** $0 < p < 1$. Let $\Omega \subset \mathbb{C}^n$ be a bounded open set. Recall that if $0 < p < 1$, we can define again the space $\mathcal{O}L^p(\Omega)$ as the set of holomorphic functions $f : \Omega \to \mathbb{C}$ such that $\int_\Omega |f(z)|^p dv(z) < +\infty$, and that with the metric

$$d_p(f,g) := \int_\Omega |f(z) - g(z)|^p dv(z), \text{ for } f, g \in \mathcal{O}L^p(\Omega),$$

$\mathcal{O}L^p(\Omega)$ becomes a complete metric space. (This follows from the fact that convergence in the space $L^p(\Omega)$ implies uniform convergence on compact subsets of $\Omega$, as we justify below.)

For a fixed $q$, with $0 < q \le 1$, we may also define the spaces

$$\bigcap_{p<q} \mathcal{O}L^p(\Omega)$$

endowed with the metric

$$\tilde{d}(f,g) := \sum_{j=1}^{\infty} \frac{1}{2^j} \frac{d_{p_j}(f,g)}{1+d_{p_j}(f,g)}, \; f,g \in \bigcap_{p<q} \mathcal{O}L^p(\Omega),$$

where $p_j$ is a sequence with $0 < p_1 < p_2 < \cdots < p_j < \cdots < q$ and $p_j \to q$ (as $j \to \infty$). Then $\bigcap_{p<q} \mathcal{O}L^p(\Omega)$ becomes a complete metric space, its topology being independant of the choice of the sequence $p_j$. In fact, a sequence $f_k$ converges to $f$, in the space $\bigcap_{p<q} \mathcal{O}L^p(\Omega)$, if and only if $d_p(f_k,f) \to 0$ for every $p < q$. In particular Baire's theorem hold in $\bigcap_{p<q} \mathcal{O}L^p(\Omega)$. Moreover we point out that the space $\bigcap_{p<q} \mathcal{O}L^p(\Omega)$, with the above topology, is also a topological vector space.

Let us recall also that if $P(a,r)$ is a polydisk, $P(a,r) = \{z \in \mathbb{C}^n : |z_j - a_j| < r_j, \; j = 1,2,...,n\}$, and $f \in \mathcal{O}(\overline{P(a,r)})$, then – by the submean value property for the function $|f|^p$ (see [13]) – we have

$$|f(a)|^p \le \frac{1}{vol(P(a,r))} \int_{P(a,r)} |f(z)|^p \, dv(z) \; (p > 0).$$

Thus if $f \in \mathcal{O}L^p(\Omega)$ and $K$ is a compact subset of $\Omega$, then choosing $\delta > 0$, sufficiently small – depending on $K$, such that

$$P_a^{\delta} := \{z \in \mathbb{C}^n : |z_j - a_j| < \delta, \; j = 1,2,...,n\} \subset\subset \Omega, \; \text{ for every } a \in K,$$

we obtain

$$|f(a)|^p \le \frac{1}{vol(P_a^{\delta})} \int_{P_a^{\delta}} |f(z)|^p \, dv(z) \le \frac{1}{vol(P_a^{\delta})} \int_{\Omega} |f(z)|^p \, dv(z), \text{ for every } a \in K.$$

This gives the well-known inequality

$$\sup_{a \in K} |f(a)|^p \le \frac{1}{vol(P_0^{\delta})} \int_{\Omega} |f(z)|^p \, dv(z).$$

In particular we see that convergence in the space $\mathcal{O}L^p(\Omega)$ implies uniform convergence on compact subsets of $\Omega$.

The following conclusions can be reached for the case «$0 < p < 1$» in the same manner as in the case «$p \ge 1$».

**Conclusions.** Theorems 4, 6, 8, 10, 12, 14, 16, and remark 17, hold also in the case $0 < q \le 1$, and remark (iii) of §5 holds for the case $0 < p < 1$, too.

**19. The spaces** $\mathrm{A}^s(\Omega)$**.** As we pointed out in §3, a totally unbounded holomorphic function in an open set $\Omega$, is singular at every point of $\partial\Omega$. On the other hand it is well-known that the converse of this is far from being correct. In fact, under some assumptions on the set $\Omega$, there are holomorphic functions in $\Omega$ which are $C^\infty$ upto the boundary of $\Omega$ and at the same time they are singular at every point of $\partial\Omega$. For deep results in this direction we refer to [8] and the bibliography given there. In this section we will use Theorem (ii) of §2 in order to give a simple proof of the fact that in some pseudoconvex open sets there exist functions in $\mathrm{A}^s(\Omega)$, $s \in \{0,1,2,...\} \cup \{\infty\}$, which do not extend holomorphically beyond any boundary point of $\Omega$. In fact we show, at the same time, that such functions form a dense and $\mathcal{G}_\delta$ set in the space $A^s(\Omega)$ (in the natural topology of this space). To make this precise, we consider, for a bounded open set $\Omega$ in $\mathbb{C}^n$ and $s \in \{0,1,2,...\}$, the set $A^s(\Omega)$ of all holomorphic functions $f$ in $\Omega$, whose derivatives

$$\frac{\partial^{|\alpha|} f}{\partial z^\alpha} = \frac{\partial^{\alpha_1 + \cdots + \alpha_n} f}{\partial z_1^{\alpha_1} \cdots \partial z_n^{\alpha_n}}$$

extend continuously to $\overline{\Omega}$, for every mult-index $\alpha = (\alpha_1,...,\alpha_n) \in \mathbb{N}^n$ with $|\alpha| = \alpha_1 + \cdots + \alpha_n \leq s$. The topology in $A^s(\Omega)$ is defined by the norm

$$\|f\|_s = \sup\left\{ \left| \frac{\partial^{|\alpha|} f}{\partial z^\alpha}(z) \right| : z \in \overline{\Omega}, |\alpha| \leq s \right\}, \ f \in A^s(\Omega),$$

and with this norm, $A^s(\Omega)$ is complete.

Similaly $A^\infty(\Omega)$ is the set of holomorphic functions $f$ in $\Omega$, whose derivatives $\partial^{|\alpha|} f / \partial z^\alpha$ extend continuously to $\overline{\Omega}$, for every mult-index $\alpha = (\alpha_1,...,\alpha_n) \in \mathbb{N}^n$. The topology in $A^\infty(\Omega)$ is defined by the metric

$$\mu(f,g) = \sum_{N=0}^\infty \frac{1}{2^N} \frac{\|f-g\|_N}{1 + \|f-g\|_N}, \ f,g \in A^\infty(\Omega),$$

and, with this metric, $A^\infty(\Omega)$ is complete. Furthermore, with the corresponding topology, $A^\infty(\Omega)$ becomes a topological vector space. Thus, in particular, if $f_k, f \in A^\infty(\Omega)$ with $\mu(f_k, f) \to 0$ ($k \to \infty$), and $\lambda_k, \lambda \in \mathbb{C}$ with $\lambda_k \to \lambda$, then $\mu(\lambda_k f_k, \lambda f) \to 0$.

The following theorem follows easily from Theorem (ii) of §2. See also [1], [4] and [8] for related results.

**20. Theorem.** *Let* $\Omega \subset \mathbb{C}^n$ *be a pseudoconvex open set such that its closure* $\overline{\Omega}$ *has a neighborhood basis of pseudoconvex open sets, and* $\text{int}(\overline{\Omega}) = \Omega$. *If* $s \in \{0,1,2,...\} \cup \{\infty\}$, *then the set* $\Xi^s(\Omega)$ *of the functions in* $A^s(\Omega)$ *which are not extentable, as holomorphic functions, beyond any point of the boundary* $\partial\Omega$, *is dense and* $\mathcal{G}_\delta$ *in the space* $A^s(\Omega)$. *In particular the conclusion holds if* $\Omega$ *is strictly pseudoconvex open set (not necessarily with smooth boundary) and* $\text{int}(\overline{\Omega}) = \Omega$.

**Proof.** We will apply Theorem (ii) of §2 with $\mathcal{V} = A^s(\Omega)$. For this purpose let us consider a pair $(B,b)$ of open balls with $b \subset\subset B \cap \Omega \neq B$. We claim that $B \cap (\mathbb{C}^n - \overline{\Omega}) \neq \varnothing$. For if $B \cap (\mathbb{C}^n - \overline{\Omega}) = \varnothing$ then $B \subseteq \overline{\Omega}$ which would imply that $B \subseteq \text{int}(\overline{\Omega})$, i.e., $B \subseteq \Omega$ (since we assume $\text{int}(\overline{\Omega}) = \Omega$), and this contradicts the fact that $B \cap \Omega \neq B$. Let $\zeta \in B \cap (\mathbb{C}^n - \overline{\Omega})$. Since we assume that $\overline{\Omega}$ has a neighborhood basis of pseudoconvex open sets, there exists a pseudoconvex open set $G$ such that $G \supset \overline{\Omega}$ and $\zeta \notin G$. Then $B \cap G \neq \varnothing$, $B \cap (\mathbb{C}^n - G) \neq \varnothing$, and $B$ is connected, and therefore $B \cap \partial G \neq \varnothing$. Let us consider a point $\sigma \in B \cap \partial G$ and a sequence $z_k$ in $B \cap G$ which converges to $\sigma$. Since $G$ is pseudoconvex, there exists a function $f$, holomorphic in $G$, such that $\sup_k |f(z_k)| = \infty$ (see [13]). Then $f \in \mathcal{V} = A^s(\Omega)$ and the restriction $f|_b$, of $f$ to $b$, has no bounded holomorphic extension to $B$. Therefore, from Theorem (ii) of §2, the set $\Xi^s(\Omega)$ is dense and $\mathcal{G}_\delta$ in the space $\mathcal{V} = A^s(\Omega)$. The last conclusion of the theorem follows from the well-known fact that the closure of a strictly pseudoconvex open set has a neighborhood basis of pseudoconvex open sets (see [5]).

**Acknowledgment.** We would like to thank M. Fragoulopoulou and A. Siskakis for helpful communications.

T. Hatziafratis, K. Kioulafa, V. Nestoridis
National and Kapodistrian University of Athens
Department of Mathematics
Panepistemiopolis
157 84 Athens
Greece

e-mail of Hatziafratis: thatziaf@math.uoa.gr
e-mail of Kioulafa: keiranna@math.uoa.gr
e-mail of Nestoridis: vnestor@math.uoa.gr